\documentclass[10pt,a4paper]{article}

\usepackage[utf8]{inputenc}
\usepackage{amsmath}
\usepackage{amsfonts}
\usepackage{amssymb}
\usepackage{graphicx}
\usepackage{color}

\usepackage{amsmath}
\usepackage{amsfonts}
\usepackage{amsthm}
\usepackage{ucs}
\usepackage[utf8]{inputenc}
\usepackage{fontenc}
\usepackage{graphicx}
\usepackage{todonotes}
\usepackage{float}
\usepackage{todonotes}

\newtheorem{theor}{Theorem}
\newtheorem{lem}[theor]{Lemma}

\newtheorem{rem}[theor]{Remark}

\newcommand{\be}{\begin{equation}}
\newcommand{\ee}{\end{equation}}

\newcommand{\bW}[0]{\mathbf{W}}
\newcommand{\bT}[0]{\mathbf{T}}

\newcommand{\mm}{\mathfrak{m}}
\newcommand{\MM}{\mathfrak{M}}

\newcommand{\Lag}{\text{Lag}}

\title{Sharpening methods for finite volume schemes}

\author{
B.~Despr\'es\thanks{Sorbonne Universit\'es, UPMC Univ Paris 06, UMR 7598, Laboratoire Jacques-Louis Lions, F-75005, Paris, France}, 
S.~Kokh\thanks{ %
  Maison de la Simulation USR 3441, Digiteo Labs, b\^at. 565, PC 190, CEA Saclay, 91191 Gif-sur-Yvette, France and 
DEN/DANS/DM2S/STMF/LMSF, CEA Saclay, 91191 Gif-sur-Yvette, France.
}
and
F.~Lagouti\`ere\thanks{Laboratoire de Math\'ematiques d'Orsay, Univ. Paris-Sud, CNRS, 
Universit\'e Paris-Saclay, 91405 Orsay, France.}
}

\begin{document}

\maketitle

\begin{abstract}
We review sharpening methods for finite volume schemes, with an
emphasis on the basic structure of sharpening methods.
It covers high order methods and non linear techniques for linear advection, Glimm's method,
anti-diffusion techniques, 
 the interaction of these techniques with the PDE structures. Additional approaches
 like level sets,   interface reconstruction and Vofire are also discussed.
 We also present the algorithmic structure of the downwind method for a simple
two components problem.  
\end{abstract}

{\bf keywords} Sharpening methods, Finite Volume schemes, anti-diffusion, interface tracking.

\section{Introduction}
The present paper deals with sharpening methods for finite volume schemes (FV)
understood as discretization strategies  for the enhancement of  sharp   profiles in numerical simulations.
We restrict the scope to finite volume schemes since they are the numerical method of choice for compressible computation fluid dynamics where the exact or approximate solutions may exhibit strong gradients that account for shocks or contact discontinuities.
%
%

We will more specifically focus on the calculations of interfaces associated with linearly degenerate fields (contact discontinuities), material contact discontinuities or free boundaries that are tracked across the computational domain. Although we shall consider numerical methods that are compliant with shock capturing, we shall not discuss the approximations of shocks in this paper.

%
%
%
Interface tracking has motivated a considerable amount of contributions since the early days of scientific computing and numerical analysis.
%
%
%
%
Therefore reviewing exhaustively all the methods that have been published to date seems quite unrealistic and we apologize in advance to the community for all the works that will not be mentioned in the sequel.
%
%
%
We propose to sketch a map of these methods by relying on mathematical and algorithmic arguments that can be used to analyze the efficiency.
We hope that this effort may also help classifying the methods that will not be discussed in this document and help understanding the sharpening mechanisms at play within the numerical schemes that are available in the literature.

%

The paper is organized as follows. 
Most of the common material (that is high order and nonlinear techniques,
the Glimm's scheme, the notion of anti-diffusion, level-sets, multidimensional issues) is presented for linear equations in section \ref{sec2}.
The introduction  of sharpening methods in nonlinear systems is evoked in section \ref{sec3}.
References are provided inside the text.





 
\section{Sharpening methods for linear equations} \label{sec2}

Sharpening methods for linear equations use two  important generic
 ideas:
the first one is to use high order schemes, and it may seem paradoxical at first sight;
the second idea is based on  compression with nonlinear techniques;  other strategies 
rely on the Glimm's scheme, on PDEs to represent the interface, or
reconstruct locally as in the volume of fluid (VOF) method.
Most of the ideas can be presented on  
the advection equation with velocity $u \in \mathbb R$, which serves as a model problem. It writes 
$$
\partial_t c(t,x) + u \partial_x c(t,x) =0, \quad x\in \mathbb R, \quad t>0, 
$$
together with a Cauchy datum $c(0,x) = c^{\rm ini}(x)$. 

\subsection{High order methods}

References to high order discretization of nonlinear equations 
are
\cite{Tor97,MR2549232,MR2337583,MR2219366}.
The fact that high order methods have the ability to sharpen discontinuities is kind of a paradox. Indeed, local Taylor expansions show bad convergence behavior for profiles involving discontinuities or strong gradients.

We give hereafter a simple  explanation of the corresponding sharpening  based on 
the theory of linear Strang's stencils. Let $\Delta t>0$ and $\Delta x$ be respectively the time and space steps. We consider a series of instants $t^n = n\Delta t$ and  the classical discretization of the real line into intervals
$[x_{j-1/2},x_{j+1/2}]$, whereby $x_j = j\Delta x$ and $x_{j+1/2} = (j+1/2)\Delta x$. We note 
$c_j^n$ an approximation of $c$ at instant $t^n$ within the cell $[x_{j-1/2},x_{j+1/2}]$ and set 
$c^n = (c_j^n)$. The initial numerical datum can be taken as $c_j^0 = c^{ini}(x_j)$ (this is especially done when dealing with smooth solutions and high order methods) or $c_j^0 = \int_{x_{j-1/2}}^{x_{j+1/2}} c^{ini}(x)\, dx/\Delta x$ (usually when dealing with non-smooth data). 
The analysis  is here limited to explicit and compact schemes with a stencil of 
$p+1$ contiguous cells. 
 In a simplified 
finite difference form 
 on a Cartesian grid, the 
family of linear schemes may read
\be \label{eq2}
c_j^{n+1}=\sum_{r=k-p}^{k} \alpha_r  c_{j+r}^n,
\qquad \alpha_r=\alpha_{r}(\nu).
\ee
The coefficients $\alpha_{r}$ are functions
of the Courant-Friedrichs-Lewy (CFL) number 
$\nu=u {\Delta t}/{\Delta x}$.
It is possible to write a scheme with order $p$ in time and space using \eqref{eq2}. 
Once $p$ has been chosen, 
$k$ determines the shift of the scheme. 
Basic examples are the well-known upwind scheme 
$c_j^{n+1}=(1-\nu)c_j^n+\nu c_{j-1}^n$, 
when $(p,k)=(1,0)$, $\alpha_{-1} = \nu$ and $\alpha_0 = 1-\nu$, 
the Lax-Wendroff scheme \cite{LaxWen60} 
$c_j^{n+1}= H^\text{LW}(c^n)_j = (1-\nu^2)c_j^n+\frac{\nu+\nu^2}2 c_{j-1}^n
+\frac{\nu^2-\nu}2 c_{j+1}^n$, 
with $(p,k)=(2,1)$, $\alpha_{-1} = (\nu^2 + \nu)/2$, $\alpha_0 = 1-\nu^2$ and 
$\alpha_1 = (\nu^2 + \nu)/2$, 
and 
the Beam-Warming scheme \cite{beam:warming}:
$c_j^{n+1}
= H^\text{BW}(c^n)_j
=\left(1-\frac32 \nu +\frac12\nu^2\right)c_j^n
+(2\nu-\nu^2) c_{j-1}^n
+\frac{\nu^2-\nu}2 c_{j-2}^n$, 
with $(p,k)=(2,0)$, $\alpha_{-2} = (\nu^2 - \nu)/2$, $\alpha_{-1} = 2\nu - \nu^2$ 
and $\alpha_0 = 1 - 3\nu/2 + \nu^2/2$.
Under the hypothesis that $\sum_{r} \alpha_r = 1$, which is a natural assumption that 
ensures the conservativity of the algorithm, 
these schemes may be rewritten also as finite volume methods in their classical form 
\be \label{eq9}
\frac{c_j^{n+1}-c_j^n}{\Delta t}+u 
\frac{c_{j+\frac12}^n-c_{j-\frac12}^n}{\Delta x}=0.
\ee
The conversion between the two forms is let to the reader because it does not 
have impact on the following discussion. 
A  third order in time and space $O3$ scheme  
$(p,k)=(3,1)$ is defined by a convex
combination \cite{despresjll}
of the Lax-Wendroff scheme and the Beam-Warming scheme:
$
c_j^{n+1} = 
(1-\alpha) H^\text{LW}(c^n)_j +\alpha H^\text{BW}(c^n)_j
$
with
$\alpha= \frac{1+\nu}3$.
The seminal works of Iserles and Strang \cite{strang,iserles1}
show that the order
in time and space, $p$, can be arbitrary large. 
Nevertheless, {\it the only pairs $(p,k)$ for which there exists schemes 
such that the $l^2$ norm is non-increasing at any iterate
for all $\nu\leq 1$ are
$p=2k+1$, $p=2k$ and $p=2k+2$.}
In the following, we call $IS$-schemes such schemes. 
\newline
The stability in $L^1$ of $IS$-schemes has been given in \cite{desstab}: 
{\it Assume moreover the order is odd, that is   $p=2k+1$. Then the scheme is   
stable in  all $L^q$:  there exists a constant $D_p >0$ such that
$
||c^n||_{L^q} \leq D_p ||c^0||_{L^q}$ $ \forall n$, $ \forall \nu \in ]0,1]$, 
$\forall c^0$ and 
$\forall q \in [1, \infty]$. 
}
\newline
Equipped with these fundamental results, a convergence result that provides a sharp
convergence estimate for an initial datum with bounded variation (BV datum) can be stated \cite{despresjll}. 
The proof is done by regularization  of the BV profile and use of the $L^1$ stability.
In this result, $c^n$ is to be understood as the constant by cell function that takes the value 
$c_j^n$ in the cell number $j$, namely $[(j-1/2)\Delta x, (j+1/2)\Delta x)$. 

\begin{theor} \label{maintheor}
Assume $c^{\rm ini}\in L^\infty\cap BV$ (in space dimension 1, this is  just the $BV$
space). 
Consider an $IS$-scheme, with $p=2k+1$ odd. 
Assume $\nu \leq 1$. 
Then
\be \label{eq:75}
||c^n-c(n\Delta t)||_{L^1}\leq C_p
|c^{\rm ini}|_{BV}
\left( \Delta x^{a }
T^{b} +\Delta x \right)
\ee
with
$
a=\frac{p}{p+1}$ and
$b=\frac{1}{p+1}$.
\end{theor}
Here, as the estimate is for non-smooth data and thus is of order less than 1, the initial numerical datum can be chosen both as point values or mean values. 

Using very high order schemes
means choosing $p$ very large.
In this case $\frac{p}{p+1}$ is very close to 1. This is 
optimal because  an error of order 1 
is what we  get by a 1 cell  translation
of the Heavyside function. In a nutshell:
very high odd order advection schemes have  nearly optimal order
of convergence in $L^1$ even
for discontinuous initial data.
It means that the very high order feature of such schemes is able to sharpen
discrete profiles with strong gradients.
Perhaps even more important for applications
is the very small 
dependence with respect to the time $T$ 
since  $\frac{1}{p+1}$ is close to zero for large $p$.
This means that the difference between the true
solution and the numerical solution
does not evolve significantly in time. That is the sharpening effect is time independent.
This theoretical behavior is the solution of the apparent paradox
explained at the beginning of the section.

Nevertheless the drawbacks of these high orders (linear) FV
methods is that they do not satisfy the maximum principle, according to a well-known theorem by Godunov. 


\subsection{Compression within a BV setting}\label{nonlinear}


It is known since \cite{harten1,harten2} and \cite{Leroux} 
that the bounded variation (BV) setting is a 
convenient framework of the construction of numerical nonlinear
FV
schemes with good sharpening properties. 
For one-dimensional problems the BV setting is strongly related 
 to the  preservation of the maximum principle \cite{roe0,sweb1,sweb2}.
In some cases the sharpening effect is so pronounced that it is called squaring. 

The general situation can be explained  a follows. Consider the FV formulation 
\be \label{eqkk1}
\frac{c_j^{n+1}-c_j^n}{\Delta t}+
u \frac{c_{j+\frac12}^n - c_{j-\frac12}^n    }{\Delta x}=0, \qquad
u>0.
\ee
Scheme~\eqref{eqkk1} 
can be recast into
\be \label{eqkk2}
c_j^{n+1} = c_j^n - \nu \left( c_{j+\frac12}^n- c_{j-\frac12}^n\right), \quad
\nu =u\frac{\Delta t}{\Delta x}.
\ee
The numerical fluxes   $c_{j+\frac12}^n$ are yet to be defined at this point of the construction. 
The design principle is to impose the maximum principle under the form
\be \label{eqkk3}
\min(c_j^n, c_{j-1}^n)\equiv m_{j-1/2}^n\leq c_j^{n+1}\leq M_{j-1/2}^n\equiv \max(c_j^n, c_{j-1}^n),
\ee
which is legitimate for advection to the right. 
 If the advection is to the left ($u<0$)  one takes 
$ \min(c_j^n, c_{j+1}^n)\equiv m_{j+1/2}^n \leq c_j^{n+1} \leq M_{j+1/2}^n \equiv \max(c_j^n,
 c_{j+1}^n)
 $.
 We consider the classical formula (see \cite{sweb1})
\begin{equation} \label{eqlw1}
 c_{j+\frac12}^n =  c_j^n +\frac12(1-\nu) 
(c_{j+1}^n -c_j^n )\varphi_{j+\frac12}^n  , \quad \forall j.
\end{equation}
where 
the correction factor  $\varphi_{j+\frac12}^n$
is a limiter or slope limiter.
It is usually defined as a function of the local slope ratio 
$$
\varphi_{j+\frac12}^n=\varphi( r_{j+\frac12}^n  ), 
\quad r_{j+\frac12}^n = \frac{c_j^n - c_{j-1}^n  }
{c_{j+1}^n -c_j^n  }.
$$
There are natural additional constraints for the definition of the slope limiter.
A first one writes 
$\varphi(1)=1$: it gives back the Lax-Wendroff flux in case $r=1$, and, generally, the second order when the datum is smooth. 
A second constraint can be 
$\varphi(r)=0$ for any $r\leq 0$: 
this is a way to enforce a local preservation of the maximum principle,
this is explained in the classical textbook \cite{Tor97}.
Another idea could be to add $\varphi(r)=r \varphi\left(\frac1r\right)$ \cite{Tor97}.
There are variants where these conditions
are relaxed, see for example \cite{dubois,Schmidtmann}.

Most of the formulas published in the literature make use of the minmod  function.
Its value is given as follows:
if   $ab\leq 0$ then   $\rm{minmod}(a,b)=0$; 
if  $a> 0$ and   $b>0$, 
then    $\rm{minmod}(a,b)=\min(a,b)$;
if 
 $a< 0$ and     $b<0$, 
then    $\rm{minmod}(a,b)=\max(a,b)$.
Then the multidimensional function 
 $\rm{minmod}: \ \mathbb{R}^p \rightarrow
\mathbb R$ is defined recursively for  $p\geq 2$ independently of the ordering by 
$$
\rm{minmod}(\mathbf{a})=
\rm{minmod}(\rm{minmod}( \mathbf{b}),c )  )
\ \  \mbox{ for }\mathbf{a}=(\mathbf{b},c)\in \mathbb R^{p},  \ \ 
\mathbf{b}\in \mathbb R^{p-1}, \ c\in \mathbb R.
$$ 
A first classical result is that if the slope limiter satisfies 
\be \label{eqlw4}
0\leq \varphi(r) \leq 
2 \rm{ minmod} \left( 1  , {r}
  \right)  
\ee
then the scheme  (\ref{eqkk1}) 
with the flux 
(\ref{eqlw1})
satisfies the maximum principle. 
%
Even if this is a very classical result, we provide a proof since
it will explain how to modify (\ref{eqlw4}) for deriving schemes with 
even stronger sharpening effect.
One has 
$$
c_j^{n+1} = c_j^n - \nu \left(c_j^n + 
 \frac12(1-\nu) 
(c_{j+1}^n -c_j^n )\varphi_{j+\frac12}^n
\right.
$$
$$
\left.
- c_{j-1}^n -  \frac12(1-\nu) 
(c_{j}^n -c_{j-1}^n )\varphi_{j-\frac12}^n \right)
$$
$$
=c_j^n
-\nu
\left(
1+\frac12(1-\nu) \left(  \frac{\varphi_{j+\frac12}^n}{r_{j+\frac12}^n}
 - \varphi_{j-\frac12}^n \right)
\right)(c_{j}^n -c_{j-1}^n ), 
$$
that is 
$
c_j^{n+1}=(1-L_j^n) c_j^n+L_j^n c_{j-1}^n, \quad
L_j^n= \nu+
\frac{  \nu (1-\nu ) }2\left(  \frac{\varphi_{j+\frac12}^n}{r_{j+\frac12}^n}
 - \varphi_{j-\frac12}^n \right)$.
The maximum principle is satisfied provided 
$0\leq L_j^n \leq 1$, that is 
$$
0\leq \nu+
\frac{  \nu (1-\nu ) }2\left(  \frac{\varphi_{j+\frac12}^n}{r_{j+\frac12}^n}
 - \varphi_{j-\frac12}^n \right) \leq 1.
$$
Assume   (\ref{eqlw4}) holds. Then   
$0\leq \varphi_{j-\frac12}^n\leq 2$
 and $1-\frac{1-\nu}2 \varphi_{j-\frac12}^n
\geq 1- (1-\nu)\geq 0$, thus $0\leq C_j^n$.
One notices that   (\ref{eqlw4}) also yields 
$0\leq \varphi_{j+\frac12}^n \leq 2 r_{j+\frac12}^n$. Therefore 
$1+\frac {1-\nu}2 \frac{\varphi_{j+\frac12}^n}{r_{j+\frac12}^n}\leq
1+(1-\nu)=2-\nu $.
Finally 
$$
\nu+
\frac{  \nu (1-\nu ) }2
 \frac{\varphi_{j+\frac12}^n}{r_{j+\frac12}^n}\leq 2\nu -\nu^2 \leq 1, \quad
\forall \nu \in [0,1], 
$$
which ends the proof.

A huge number of formulas has been proposed in the literature.
We just review the most usual ones.
The Minmod flux writes 
\be \label{eqminmod}
\varphi(r)= \rm{minmod}(1,r).
\ee
The Superbee flux writes
\be \label{eqsuperbee}
\varphi(r)= \max(0,\min(1, 2r   ), \min(2, r)).
\ee

\begin{rem}[Squaring/sharpening behavior of  Superbee] \label{rem2}
The notion of sharpening is not present at this stage of the discussion.
It is introduced by noticing that the SuperBee limiter is squaring.
This has been reported in the literature in \cite{Tor97} and many other texts.
Squaring means that if an initial smooth profile is chosen, for example in the form of a Gaussian,
then the numerical solution has the tendency to converge to a mass
preserving square profile for $t\rightarrow \infty$. This behavior necessarily increases the $L^2$ norm of the profile. 
\end{rem}

 Squaring  is usually  considered as a
 consequence of the strong nonlinearity of SuperBee.
 Even if it is a well documented behavior, we know of no definitive  proof.
 But on the contrary,  it is easy to understand that the minmod  limiter 
cannot sharpen.
 To this end we consider the semi-discrete (that is continuous in time)
version of the scheme 
\begin{equation} \label{eq:mm}
\frac {\text{d}}{\text{d}t} c_j(t) + u \frac{c_{j+\frac12}(t) - c_{j-\frac12}(t)    }{\Delta x}=0, \quad j\in \mathbb Z.
\end{equation}
Since $ \Delta t$  vanishes, the flux  (\ref{eqlw1}) 
 is simplified taking  $\nu=0$. 

\begin{lem}
The semi-discrete scheme (\ref{eq:mm}) with the flux $ c_{j+\frac12} =  c_j +\frac12 
(c_{j+1} -c_j )\varphi_{j+\frac12}$    and the Minmod limiter (\ref{eqminmod})
satisfies the a priori estimate
\begin{equation} \label{eq:mm2}
\frac d {dt } \left( \sum_{j\in \mathbb Z} |c_j(t)|^2 \right)   \leq 0.
\end{equation}
So, as a corollary of remark \ref{rem2},  this scheme  cannot sharpen. 
\end{lem}

The same  property holds for similar schemes with a limiter $0\leq \varphi(r)\leq 1$ for all $r$.
The proof proceeds as follows. One has 
$$
\frac{\Delta x}{2} \frac {\text{d}}{\text{d}t}\left( \sum_{j\in \mathbb Z} |c_j|^2 \right)= \Delta x \sum_j c_j \frac {\text{d}}{\text{d}t} c_j=
- u  \sum_j c_j \left(  c_{j+\frac12} - c_{j-\frac12}  \right)
$$
$$
= - u  \sum_j c_j \left(  c_j - c_{j-1}  \right)- u \sum_j c_j \left(  c_{j+\frac12} - c_j  \right)+
u\sum_j c_j \left(  c_{j-\frac12} - c_{j-1} \right).
$$
It is easy to check the identities 
$$
\left\{
\begin{array}{lll}
\sum_j c_j \left(  c_j - c_{j-1}  \right) &= \frac12 \sum_j |c_j-c_{j-1}|^2, \\
\sum_j c_j \left(  c_{j+\frac12} - c_j  \right)&= \frac12 \sum_j c_j (c_{j+1} -c_j )\varphi_{j+\frac12}, \\
&= \frac12 \sum_j c_{j-1} (c_{j} -c_{j-1} )\varphi_{j-\frac12}, \\
\sum_j c_j \left(  c_{j-\frac12} - c_{j-1}  \right)&= \frac12 \sum_j c_j (c_{j} -c_{j-1} )\varphi_{j-\frac12} .
\end{array}
\right.
$$
Therefore by summation and rearrangements 
$$
\frac{\Delta x}{2} 
\frac {\text{d}}{\text{d}t}\left( \sum_{j\in \mathbb Z} |c_j|^2 \right)=- \frac u 2  \sum_j |c_j-c_{j-1}|^2 \left( 1- \varphi_{j-\frac12} \right)
\leq 0, 
$$
which shows that the $L^2$ norm decreases. It makes squaring impossible.
The proof is ended.
It can be generalized to the fully discrete scheme with the same conclusion. 
A corollary is as follows.

\begin{lem} [Necessary condition for sharpening]
A slope limiter that sharpens is necessarily such that 
$\varphi(r)>1$ for some $r \in \mathbb{R}$. This condition is satisfied by  the SuperBee formula
(\ref{eqsuperbee}), for which
$\lim_{r\rightarrow \infty}\varphi(r)=2$.
\end{lem}

\subsection{Inequality and anti-diffusion}\label{section : inequality and anti-diffusion}

This sharpening strategy is more radical.
It is naturally introduced in the context of BV schemes \cite{DesLag02},
see also  \cite{Tor97},
and has been adapted to ENO techniques \cite{MR2481112} in \cite{cws1,cws2}.
We refer to \cite{shyue-xiao,MR2877636,MR2787572,KoLag2010,MR2481112,jaouen-lagoutiere,cws1,cws2,billaud-kokh,goudon-lagoutiere-tine} 
for the use of such methods for different problems.

We shall note $ m_{j+1/2}= \min(c_j^n,c_{j+1}^n) $,  $ M_{j+1/2}= \max(c_j^n,c_{j+1}^n) $
\begin{equation}
 \lambda_{j+1/2} =
 \frac{\Delta x}{u \Delta t}
(
c_j^n
-
M_{j-1/2} 
)
 +
M_{j-1/2}^n
\mbox{ and }
\Lambda_{j+1/2}
=
\frac{\Delta x}{u \Delta t}
(
c_j^n
-
m_{j-1/2} 
)
 +
m_{j-1/2}^n
  .
\label{eq: constant velocity def of a and A}
\end{equation}
We observe that  $\lambda_{j+1/2} \leq  \Lambda_{j+1/2} $ if the CFL condition $u\Delta t \leq \Delta x$ is satisfied.
A basic property writes as follows.

\begin{lem}
Under CFL,  the upwind flux choice  $c^n_{j+1/2} = c^n_j$ belongs 
 to the interval $[ \lambda_{j+1/2}, \Lambda_{j+1/2}] \cap [m_{j+1/2},M_{j+1/2}]$, which ultimately ensures
the maximum principle \cite{DesLag02,cws1,cws2}. 
\end{lem}

In this context, one  introduces compression, or sharpening, or anti-diffusion, 
by using the most extreme formulated choice.  Let $\omega_{j+1/2}$ and $\Omega_{j+1/2}$ such that 
$[\omega_{j+1/2},\Omega_{j+1/2}] = [ \lambda_{j+1/2}, \Lambda_{j+1/2}] \cap [m_{j+1/2},M_{j+1/2}]$. One obtains
\begin{equation*}
c_{j+1/2}^n=
 \begin{cases}
  \Omega_{j+1/2},& \text{if $ \Omega_{j+1/2} \leq c_{i+1}^n $,}
\\
  c_{j+1}^n,& \text{if $ \omega_{j+1/2} \leq c_{i+1}^n \leq \Omega_{j+1/2} $,}
\\
  \omega_{j+1/2},& \text{if $ c_{i+1}^n \leq \omega_{j+1/2}$.}
\end{cases}
\end{equation*}
An equivalent definition (still for the case $u > 0$) is given in the following lemma. 

\begin{lem}\label{equivalenceUltrabee}
The limited downwind flux defined above is equivalent to the so-called 
UltraBee flux limiter flux  (see \cite{Tor97}) defined as 
$$
c_{{j}}^{{n+1}} =
c_{{j}}^n-\nu(c_{{j}}^n-c_{j-1}^n) - 
\frac{\nu(1-\nu)}{2}
(\varphi_{{j+1/2}}^n(c_{{j+1}}^n-c_{{j}}^n) -
\varphi_{{j-1/2}}^n(c_{{j}}^n-c_{j-1}^n))   
$$
with
$
\varphi_{{j+1/2}}^n=\varphi(r_{{j+1/2}}^n,\nu)$ 
and 
$
\varphi(r,\nu) =
\mbox{minmod}(\frac{2r}{\nu},\frac{2}{1-\nu}) 
$.
\end{lem}

The limiter is now function of the slope $r$ and of the Courant number $\nu$. The scheme is called {\em limited downwind} in the following. 

\begin{lem}
This limited downwind scheme is exact for step initial conditions \cite{DesLag02}. 
\end{lem}

Confirmation is by starting from an initial data which is not a step function, but a (discretized) smooth function.
One observes (under a surprising technical condition CFL$\neq 1/2$) that the smooth profile is replaced
a step function close by a step function with an  approximation error is  $O(\Delta x)$. After that first stage the step function
is perfectly transported. 
So in some sense the UltraBee limiter is a perfect sharpener.
The sharpening effect  is so pronounced that it may resemble an instability, but it is not.

This technique was incorporated in FV algorithms for the simulation of two-component fluid flows, for the mass fraction, volume fraction, or color function of components, in, e.g., \cite{DesLag2007}, \cite{KoLag2010}, and extended to {\em multi}-component in \cite{jaouen-lagoutiere} and \cite{billaud-kokh}. 

This was also modified to apply to non-linear discontinuities such as classical shocks, in \cite{aguillon-chalons}, and non-classical shocks in the scalar context, \cite{boutin-chalons-lagoutiere-lefloch}, and in the context of systems in \cite{aguillon}. 

\subsection{Glimm's method}\label{section: Glimm}


At this stage of the discussion the problem is the following: either one accepts to violate the maximum principle although this can be very critical, for example when the transported unknown is the mass or volume fraction of a fluid in a multi-component flow, or one has to use linear first order or a non-linear scheme (see Section~\ref{nonlinear}). 
Yet there exists an alternative, that was first proposed by Glimm, in \cite{glimm} for theoretical analysis purposes. This method avoids the numerical diffusion of first order stable schemes because it does not involve any "projection" on the mesh, and it does not create new values of the solution a each time step in the case of linear transport. 

To describe it briefly, let us consider once again the upwind scheme 
written as $
c_j^{n+1} = (1 - \nu)c_j^n + \nu c_{j-1}^n$ with $\nu$ the CFL number.
The smearing of the profiles comes from the (strictly) convex combination that appears in the formula. This scheme can be interpreted as a two-step scheme: exact transport of the profile for a time $\Delta t$, and then projection on the mesh (the upwind scheme is the Godunov scheme). 
Glimm proposes to avoid the projection by taking one of the two values that are present in cell $j$ after one time step: $c_{j-1}^n$ or $c_j^n$. The choice is performed randomly: $c_{j-1}^n$ is chosen with probability $\nu$, and $c_j^n$ is chosen with probability $1 - \nu$. 
This interpretation is correct since  $0< 1 \leq \nu$ (resp. $
0\leq 1-\nu<1$) under CFL. In the more general context of nonlinear problems, the algorithm is based on the resolution of the Riemann problems at each interface and on the choice of a random variable $\delta^n$(different from one time step to the other), chosen according to the uniform law between $0$ and $\Delta x$. Then the updated value of the unknown in the cell $j$ is defined by taking the value of the solution of the Riemann problem at time $\Delta t$ at position $(j-1/2)\Delta x + \delta^n$. This was shown by Glimm 
to converge, with probability 1, and it is clear that it does not {\em smear} profiles, at least when the profile is a step and in the linear context. Let us note that this random procedure has the drawback that the scheme is non-conservative,  however this does not prevent the scheme to converge to the entropy solution for nonlinear problems. Note also that the randomness is not mandatory: the only property that is required for $(\delta^n)_n$ is that it is an equidistributed (with low discrepancy) sequence. The Van der Corput sequence, which is such a deterministic sequence, is shown to give qualitatively very good (better than a random sequence) results in \cite{colella}. 

One can notice that, for the linear transport equation~\eqref{eq2}, the upwind scheme is the {\em expectation} of Glimm's scheme. This observation was used to prove error estimates for the upwind scheme on general meshes, using central-limit type estimates, in \cite{delarue-lagoutiere}. 

In space dimension 1 and in the context of linearly degenerate fields (which correspond to material discontinuities)  an FV algorithm based on a Lagrange-remap (formulated as Lagrange-transport) strategy with a random sampling technique for the transport part, for the simulation of two-component compressible fluid flows is derived in \cite{chalons-coquel}, with very good efficiency. See also \cite{chalons2007,ChalonsGoatin2008,Bachmann2013275}. 

In the more particular context of {\em non-linear} material discontinuities that are present in some viscous-dispersive limits of systems with fields that are neither genuinely non-linear nor linearly degenerate, with so-called {\em non-classical shocks} (see \cite{lefloch}), the random choice method was shown in \cite{chalons-lefloch} to give very good (and convergent) results, which is very difficult in this context. 

The tentation to use such a scheme in higher dimension is great, but it is known since Chorin (\cite{chorin}) that it is not satisfactory 
for genuinely nonlinear conservation laws. Colella proposed in \cite{colella} a modification of the random choice algorithm, that involves the Godunov method and that seem to be convergent. Unfortunately this modification of Glimm’s algorithm does not allow to preserve sharp fronts. 
Nevertheless, for linear or linearly degenerate fields, this random choice procedure shows great efficiency, at least on Cartesian grids. This has been investigated and analyzed, for example in \cite{helluy-jung1} and \cite{helluy-jung2}. 

\subsection{PDE models and sharpening methods}

Level sets methods are discussed  in \cite{osher,osher-sethian,sethian}. This is a very popular set
of numerical methods for interface modeling that has been applied 
to many problems. In the present context, the idea is to rely on a partial differential equation to transport
a color function (our definition of a color function $f$ is that it takes value in $[0,1]$, so that:  if  a point $x$ is  such that
$f(x)=0$ has color equal to 0; if $f(x)=1$ then $x$ has a color equal  to 1; and finally $0<f(x)<1$ corresponds to intermediate colors).
No colors below 0 and above 1 are considered in this presentation, but it is not mandatory.
A typical elementary question with the level set approach is about the influence of the numerical parameters
on the level set.
In certain cases the answer is that the method can be insensitive to this parameters. 

To understand this property we consider the simplest color function  at initial time
$c^{\rm ini}(x)=H(x)$, that is $c^{\rm ini}(x)=0$ for $x<0$ and $c_0(x)=1$ for $x>1$.
Instead of manipulating  the upwind first order scheme, we use its modified equation (that is to say, the PDE it is consistent with at the { \em second} order in time and space). We thus consider the function $c_\mu$ solution of the advection  equation with viscosity 
$$
\partial_t c_\mu+ u \partial_x c_\mu =\mu  \partial_{xx} c_\mu , \quad
\mu= \frac{\Delta x}{2} (1-\nu), 
$$ 
where $0\leq \nu \leq 1$ is the CFL number.  The modified equation is a second order approximation of the upwind scheme.
The interface is recovered at any time $t$ by as the $1/2$ level set 
$
\Gamma_\mu(t)=
x \mbox{ such that  }c_\mu(x,t)=1/2$. It is easy to prove that $x$ exists and is unique 
for $t>0$ and $0\leq \nu<1$:  this is a consequence
of well known integral representation formula detailed below.
One has more.

\begin{lem}
For all $t>0$ and $0\leq \nu<1$, 
the $1/2$ level set is exact: that is $\Gamma_\mu(t)=ut$.
\end{lem}
 One has with the fundamental solution of the heat equation 
 $$
 c_\mu(x,t)= \frac1{\sqrt{4\pi \mu t}}
 \int_{\mathbb R} \exp\left( -(x-y-ut)^2/(4\mu t)  \right) H(y)dy.
 $$
 So
 $$
 c_\mu(ut,t)=\frac1{\sqrt{4\pi \mu t}}
 \int_0^\infty \exp\left( -y^2/(4\mu t)  \right)dy=
 \frac1{\sqrt{2\pi}}
 \int_0^\infty \exp\left( -y^2/2  \right)dy=\frac12.
 $$

Even if this argument is very elementary, it explains that level set methods have the ability
to predict the interface with great accuracy, even if the underlying
scheme for the transport of the color function is low order.
In the context of this review paper, it is perfect sharpening.



\subsection{Nature of the grid/mesh}

The discussion so far was restricted 
to one dimensional grids.
The extension of the previous FV algorithms
to general multidimensional grids  poses two fundamental difficulties.

A first one is that sharpening techniques are highly nonlinear methods.
A good sharpening technique is in practice equipped with
a method which controls the oscillations due to strong nonlinear interactions.
In dimension $d=1$, this principle is mostly based on the BV setting.
The issue is that this bound on the total variation is lost in dimension $D=2$ 
and greater. This has been proved
in a famous article \cite{good-lev} on a Cartesian grid. This unfortunate situation has the consequence
that the preservation of the maximum principle does not yield a control of some special oscillations
which develop mostly tangentially to the isolines of the exact profiles: an important reference in this direction
is the series  \cite{kuz-tur,kuz-tur2,kuz-tur3}.
See also \cite{vofire}.

A confirmation of  this behavior is the 2D algorithm in \cite{DesLag02}.
It is shown that the extension of the Ultra-Bee scheme with directional splitting is exact
for squares. But unfortunately this algorithm is not equipped with a control of 2D  variations.
It can be interpreted as a distant consequence of the 
\cite{good-lev} theorem. In consequence 
this algorithm  is useless for calculations of  profiles with values which are not exactly 0 and or 1. 
Even if the initial data is an indicatrix function, its boundary is not necessarily  a 2D step function:
in this situation one observes oscillations at the boundary between 0 and 1: these oscillations are perfectly bounded
in $L^\infty$ norm because directional splitting preserves the maximum principle;
but they are not bounded in the BV semi-norm because the BV semi-norm  is a global 
quantity destroyed by directional splitting.
An attempt is been made in \cite{lvd} to overcome this failure, but the numerical results are deceptive (not published),
probably due to the curse explained by the \cite{good-lev} theorem.

\subsection{Interface reconstruction and VOF}

The simple line interface calculation (SLIC) \cite{Slic}
is an extremely popular method that presents a nearly all purpose
methodology for FV interface sharpening.
The design principle of SLIC is to reconstruct parallel and/or anti parallel perfect 
interfaces in Cartesian cells from the knowledge of volume fractions. In dimension 1, for a step initial conidtion, it is equivalent to the limited downwind scheme (that can be seen as a reconstruction algorithm, where the reconstructed solution is a step function in every cell). 
Even if it is an extremely simple method, the results are quite good 
when comparing with the implementation cost and run time.
This is probably the reason why it is still a reference.
With respect to SLIC, the volume of fluid (VOF)  \cite{vof} has the huge advantage
to reconstruct  interface with any direction. Even without discussing the simplicity
of the method, it is clear that this information is a kind of first  order interface reconstruction 
while SLIC can be considered as a zeroth order interface reconstruction.
Another feature of VOF  is that the normal direction of the interface 
is computed from the discrete gradient of some volume fractions. 
It is  possible to optimize the performance of VOF by changing the parameters
of the discrete gradient operator and of the method  used to evolve the volume fractions.

It must be noted that  SLIC is not   PDE based and VOF is only partially PDE based. In consequence 
it is not really possible to perform a convergence analysis of the algorithms, but only on parts
of them.
The Youngs algorithm \cite{youngs} has a similar nature. 

\subsection{Vofire}

\begin{figure}[H]
\begin{center}
\scalebox{0.5}{\begin{tabular}{cc}
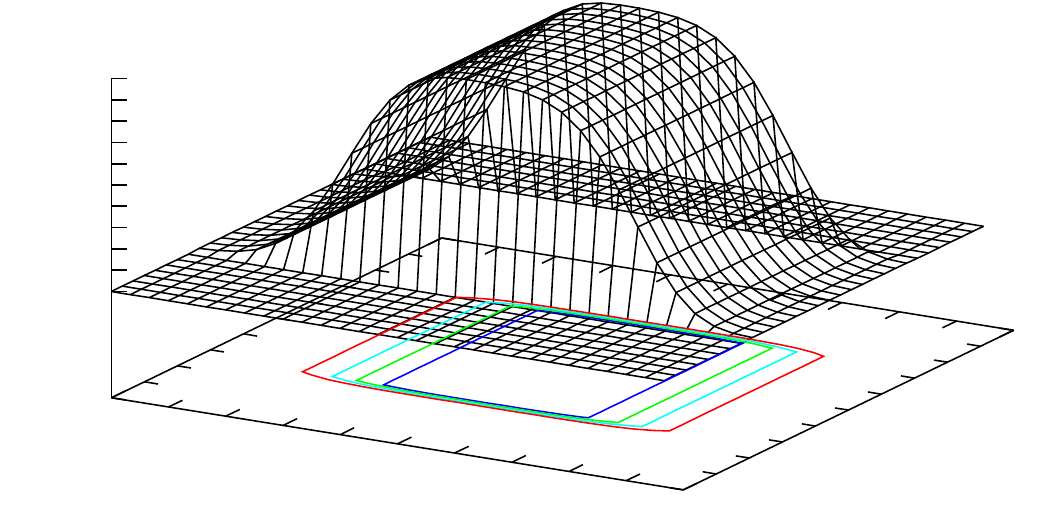 &
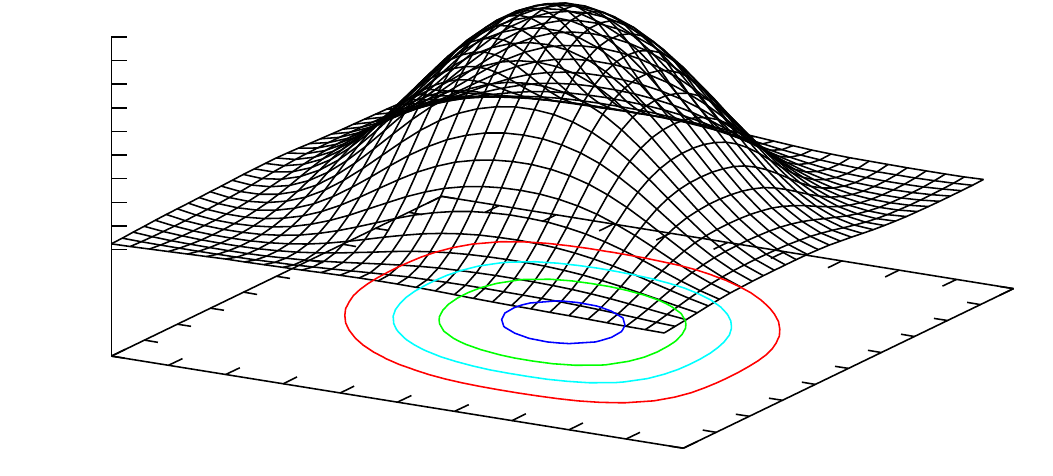
\end{tabular}}
\caption{Upwind scheme. The initial condition
is the characteristic function of a square.
Final time $t=1$. Periodic boundary conditions. 
On the left: the velocity $\mathbf{u} =(1,0)^T$ is aligned with the
mesh; the result displays only longitudinal diffusion. On the right:
the velocity $\mathbf{u} =(1,1)^T$ is not aligned with the mesh. The
consequence is that there is both longitudinal and transverse diffusion.} 
\label{longi}
\end{center}
\end{figure}

We give some details of the Vofire method, which is a multidimensional
nonlinear  FV scheme.
The geometrical idea  relies on the following observation: in dimension greater than 2, 
the numerical diffusion can be decomposed into two different diffusions: the longitudinal diffusion, 
along the velocity field, which is typically one-dimensional, and the transverse diffusion, which is really 
due to the fact that the mesh is multi-dimensional. 
This distinction between
the two phenomena could appear arbitrary, but is in accordance with
basic numerical tests. Consider for example an initial condition 
which is the characteristic function of the square
$]0.25,0.75 [\times ]0.25,0.75 [$.
This profile is advected with the upwind scheme. The velocity direction
$\mathbf{u}$ has a great influence on the result. It is
illustrated on figure \ref{longi}.

We here propose to restrict to triangular meshes, on which it is simpler to expose the Vofire technique. 
Thus we consider the following type of mesh structure: 
\begin{figure}[H]
\begin{center}
\scalebox{0.6}{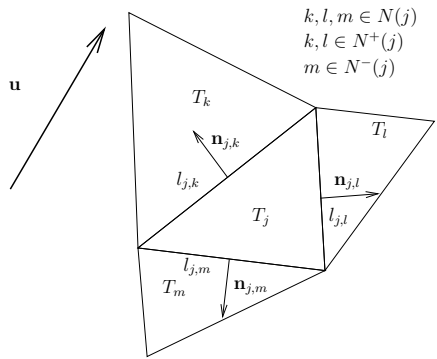}
\caption{Mesh and notations (for the sake of brievity, these notations will not be explained further in the text, as they are very classical). }
\label{figmaill}
\end{center}
\end{figure}

The idea to ''fight'' against these two diffusion phenomena is to use, as for the limited downwind seen as a 
(discontinuous) reconstruction scheme, a reconstruction procedure 
in every cell. This reconstruction will be here two-fold: it will consist in a first reconstruction that will 
be transverse, and in a second one that will be done along the velocity field.
The velocity field  $\bf u$ is constant for simplicity,  but this 
assumption can  be removed. 
Note that the transverse diffusion actually depends more on the shape of the cells than on the velocity, so that 
this assumption of constant velocity is not so much restrictive here. For a given celle 
$T_j$, $N^-(j)$ denotes the set of cells that are adjacent to $T_j$ and {\em upwind}: 
$N^-(j) = \{ T_m \mbox{ such that } \overline{T_j} \cap \overline{T_m} \mbox{ is of non-zero one dimensional Lebesgue measure and } u \cdot n_{j,m} < 0\}$, and $N^+(j)$ denotes the set of downwind cells to $T_j$. 
As in dimension 1, the fundamental requirement of the scheme is that it satisfies an upwind maximum principle: 
$$
\min\left( c_j^n, \min_{k \in N^-(j)}c_k^n \right)
\leq c_j^{n+1} \leq 
\max\left( c_j^n, \max_{k \in N^-(j)}c_k^n \right) \quad \mbox{for any } j. 
$$
The most important part of the procedure, regarding the multidimensional properties of the scheme, is the first one, 
that concerns the transverse reconstruction. As we will see, after this reconstruction, the algorithm will 
be one-dimensional, and one-dimensional techniques (such as the limited downwind scheme for instance) will  be 
applied. 

\begin{figure}[H]
\begin{center}
\scalebox{0.6}{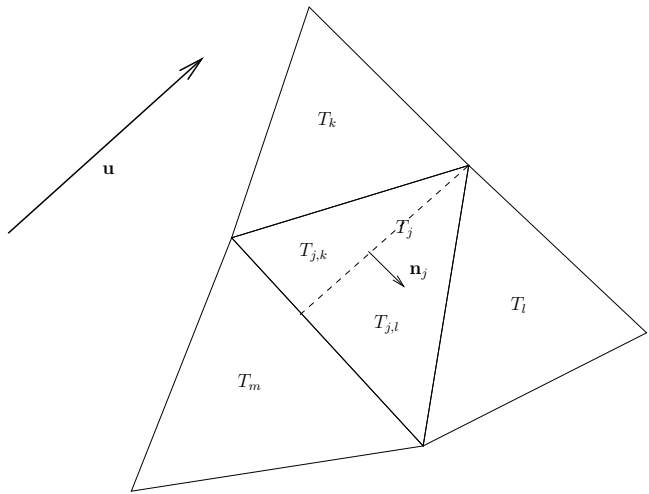}
\caption{Transverse reconstruction. }
\label{decoup}
\end{center}
\end{figure}

Recall that, for expository purposes, the mesh is assumed to be made with triangles, 
in dimension 2. 
The transverse
reconstruction consists in breaking a cell in two parts by a segment
parallel to the velocity, and modifying the value of the unknown
 in each
of these two sub-cells. Each triangle $T_j$ has at least one downwind
neighbor and at most two. If it has only
one downwind neighbor, we do not perform the transverse reconstruction
(we do not cut the cell). This can be explained by the fact that when 
there is only one downwind neighbor, the ''information'' contained in 
the cell is not spread transversally by any scheme (with small stencil). 
Let us thus assume that $T_j$ has two
downwind neighbors, $T_k$ and $T_l$. It has then one upwind
neighbor, $T_m$. We consider the intersection point of the two
edges relative to the downwind neighbors and cut $T_j$ along the line
passing on this intersection point and parallel to $\mathbf u$. The
two sub-cells are denoted $T_{j,k}$ and $T_{j,l}$: $T_{j,k}$
has $T_k$ as (unique) downwind neighbor, and $T_{j,l}$ has $T_l$ as
(unique) downwind neighbor. This partitioning is illustrated on
figure~\ref{decoup}. 
The essential property of this cutting is that now every subcell has 
exactly one downwind and one upwind neighbor, as will be used below: 
this is due to the fact that the new normal vector $\mathbf{n}_j$ 
(see Figure \ref{decoup}) is orthogonal to the velocity, so that there will 
be no flux through the new interface. 
We use symbols $s_{j,k}$ and $s_{j,l}$ to denote the areas 
 of sub-cells $T_{j,k}$
and $T_{j,l}$ respectively. Clearly $s_{j,k} + s_{j,l} = s_j$ and
$s_{j,k} > 0$ and $s_{j,l} > 0$. 
The aim is to define a reconstructed value $c_{j,k}^R$ in $T_{j,k}$
and a reconstructed value $c_{j,l}^R$ in $T_{j,l}$.
We impose 
\begin{equation} \label{volfi6}
s_{j,k} c_{j,k}^R + s_{j,l} c_{j,l}^R = s_j
c_j^n
\end{equation}
to guarantee the local conservativity.
Let us write
\begin{equation} \label{llam}
\left\{
\begin{array}{ll}
c_{j,k}^R= c_j^n + \lambda_{j,k} \left( c_k^n -c_j^n \right), & 0
\leq \lambda_{j,k} \leq 1, \\ 
c_{j,l}^R= c_j^n + \lambda_{j,l} \left( c_l^n -c_j^n \right), & 0
\leq \lambda_{j,l} \leq 1,
\end{array}
\right.
\end{equation}
which means that $c_{j,k}^R$ and $c_{j,l}^R$ must satisfy a 
condition of local consistency. We introduce the idea of 
anti-dissipative schemes which will serve to find a unique value
of $ \lambda_{j,k}$ and $ \lambda_{j,l}$. 

\begin{rem}\label{rem1}
Our objective is to choose $ \lambda_{j,k}$ and $ \lambda_{j,l}$ in
order to obtain an anti-dissipative scheme, with a very low level
of numerical diffusion.
This is the reason why we seek
the largest possible $\lambda_{j,k}$ and $\lambda_{j,l}$ in the
interval $[0,1]$. This is the same principle as in section \ref{section : inequality and anti-diffusion}.
But of course we cannot take $ \lambda_{j,k}= \lambda_{j,l}=1$
directly because we ask the reconstruction to be conservative: 
see Equation (\ref{volfi6}).
\end{rem}

So our goal is to have the largest $ \lambda_{j,k}$ and $
\lambda_{j,l}$ in the interval $[0,1]$, but still satisfying the
maximum principle. Equation (\ref{volfi6}) can be rewritten as
$s_{j,k} c_{j,k}^R + s_{j,l} c_{j,l}^R - s_j
c_j^n=
s_{j,k} \left( c_{j,k}^R -c_j^n \right) + s_{j,l} \left(
c_{j,l}^R - c_j^n\right)=0$, that is 
$$
 \left[ s_{j,k} \left( c_k^n -c_j^n \right)\right] \lambda_{j,k}+
 \left[s_{j,l}\left( c_l^n -c_j^n \right)\right]\lambda_{j,l} =0.
$$ 
As we will see, either the data $c_j^n$ is a local transverse maximum 
or minimum and then (\ref{volfi6}) implies $\lambda_{j,l} = \lambda_{j,k} = 0$ 
(it means  there is actually no reconstruction), or the data is transversally 
monotone, and then at least $\lambda_{j,l} = 1$ or $\lambda_{j,k} = 1$. 
The solution is computed as follows.

\begin{description}
\item [1) If $s_{j,k}\left( c_k^n -c_j^n \right)
s_{j,l}\left( c_l^n -c_j^n \right)\geq 0$, ]
$c_j^n$ is a local extremum in the transverse direction. Then
 we
do not reconstruct, which means $ \lambda_{j,k}= \lambda_{j,l}=0$ and 
\begin{equation} \label{eqxs1}
c_{j,l}^R = c_{j,k}^R = c_j^n.
\end{equation}
\item[2) If $\displaystyle -\frac{s_{j,k}\left( c_k^n -c_j^n \right)}
{s_{j,l}\left( c_l^n -c_j^n \right)}>1$,] the solution is obtained by
taking $ \lambda_{j,l}=1$, 
\begin{equation} \label{eqxs2}
 c_{j,l}^R = c_l^n, \qquad
 c_{j,k}^R = 
c_j^n - \frac{s_{j,l}}{s_{j,k}}\left(c_l^n -c_j^n \right)=
(s_j c_j^n - s_{j,l}
c_l^n)/s_{j,k}. 
\end{equation}
\item[3) If $\displaystyle -\frac{s_{j,k}\left( c_k^n -c_j^n \right)}
{s_{j,l}\left( c_l^n -c_j^n \right)}<1$,] the solution is obtained by
taking $ \lambda_{j,k}=1$, 
\begin{equation} \label{eqxs3}
 c_{j,k}^R = c_k^n, \qquad
 c_{j,l}^R = c_j^n - \frac{s_{j,k}}{s_{j,l}}\left(c_k^n -c_j^n
\right) = (s_j c_j^n - s_{j,k} c_k^n)/s_{j,l}.
\end{equation}
\end{description}

As the situation is now one-dimensional for each subcell $T_{jk}$ and $T_{jl}$, 
one can analyze the scheme where the first stage is this reconstruction
followed by a second stage which 
is  the upwind scheme. 
It is obvious that this scheme will provide the maximum principle, as the 
reconstruction does. What is not so obvious  is that the CFL 
stability condition for the upwind scheme on this new (and finer) mesh {\em is the 
same as for the initial mesh}. 
%
A simple proof is as follows.
\newline
{\bf Proof.}
The reconstructed
quantities (\ref{llam}) respect the maximum principle. 
By construction  the scheme is equal to a two-steps algorithm:
first step, use the upwind scheme
 for a mesh which is locally cut in smaller cells, as it is
described in figure \ref{decoup}, and with cell quantities equal to the 
reconstructed
quantities;
second step, project onto the original coarse mesh. Therefore it is sufficient
to check that the  CFL condition is the same for
the original mesh (4 cells in figure \ref{decoup}) and for the new
mesh (5 cells in figure \ref{decoup}).

Since $\mathbf u$ is constant, then 
$
\sum_{k \in N^+(j)}
l_{j,k} (\mathbf{u}^T \mathbf n_{j,k}) = -\sum_{k \in N^-(j)}
l_{j,k} (\mathbf{u}^T \mathbf n_{j,k}).
$
The standard CFL condition  for the upwind
scheme for the cell $T_j$ thus has the form 
$
\frac{\Delta t}{s_j} \sum_{k \in N^+(j)} l_{j,k} 
(\mathbf{u}^T\mathbf n_{j,k}) \leq 1$,
that is 
\begin{equation}\label{cfl}
\frac{\Delta t}{s_j} \left( l_{j,k}  (\mathbf{u}^T\mathbf n_{j,k})
+ l_{j,l} (\mathbf{u}^T\mathbf n_{j,l})\right) \leq 1. 
\end{equation}
The CFL condition for the sub-cells $T_{j,k}$ and $T_{j,l}$ are 
respectively
\begin{equation}\label{sscfl}
\frac{\Delta t}{s_{j,k}} l_{j,k} (\mathbf{u}^T\mathbf n_{j,k})
\leq 1 \quad \mbox{and} \quad \frac{\Delta t}{s_{j,l}}  l_{j,l}
(\mathbf{u}^T\mathbf n_{j,l}) \leq 1. 
\end{equation}
Let $l_j = \mbox{length} \left(\overline{T_{j,k}} \cap
\overline{T_{j,l}}\right)$ be the length of the segment separating
$T_{j,k}$ and $T_{j,l}$. One has 
$
s_{j,k} = \frac{l_j}{2 \left| \mathbf u
\right| } 
l_{j,k} (\mathbf{u}^T\mathbf n_{j,k}) $
 and 
$s_{j,l} = \frac{l_j }{2 \left| \mathbf u \right|}
 l_{j,l} (\mathbf{u}^T\mathbf n_{j,l})
$
 and 
$
s_j = s_{j,k} + s_{j,l} = \frac{l_j }{2\left| \mathbf u \right|} 
 \left(l_{j,k} (\mathbf{u}^T\mathbf n_{j,k}) + l_{j,l} (\mathbf{u}^T\mathbf
n_{j,l})\right)
$. 
The two inequalities of (\ref{sscfl}) and inequality (\ref{cfl}) thus
rewrite 
$\left| \mathbf u \right| \frac{2\Delta t}{l_j} \leq 1$. 
So they
are equivalent and  the proof is ended. 

Some modifications and improvements of the Vofire technique have been proposed in \cite{tran1, tran2} for example. 

\section{Coupling with hyperbolic  nonlinear equations} \label{sec3}

An issue is to use  the previous techniques
in complex computational fluid dynamics FV codes
with a respect of the key properties
necessary for  a correct simulation.
Of course 
 the notion of a correct simulation and the identification of the key ingredients
of a correct coupling are far to be evident. 
We restrict the discussion to hyperbolic models for compressible fluids for which 
conservative issues are critical. Indeed the nonlinearity of the equations  induce 
discontinuous solutions such as shocks and contact discontinuities: it is well 
accepted that the violation of the conservation is only at the cost of a strong deviation
with respect
to the solution of the Riemann problem (see \cite{hou-lefloch} for a justification). 
These questions are fiercely debated when 
dealing with multi-fluid models. 

There are common guidelines for incorporating sharpening techniques into discretization strategies of complex models. Usually one singles out transport effects in the system and update a set of key fluid parameters thanks to a numerical scheme that transports discontinuities as sharply as possible. A delicate matter is generally to preserve good stability and consistency properties of the overall numerical scheme.

\subsection{An example of discretization for compressible flows with two components separated by a sharp interface}
\label{section: system with dalton law}

For the sake of illustrating these ideas, we consider, in space dimension 1, a simple model of compressible flows that involves two perfect gases that was 
studied in \cite{abgrall2,larrouturou1989}. The specific heat at constant volume and ratio of specific heat of the fluid $k=1,2$ are respectively ${c_{v}}_k > 0$ and $\gamma_k > 1$. The density of the two-phase medium is noted $\rho$ and $Y_1 = Y$ (resp. $Y_2 = 1-Y$) is the mass fraction of the fluid $k=1$ (resp. $k=2$). We suppose that there is a thermal equilibrium between the gases and that the pressure $P$ verifies Dalton's law, then we have
\begin{equation}
P = 
\frac{\sum_{k=1,2} Y_k (\gamma_k -1) {c_v}_k}{\sum_{k=1,2} Y_k  {c_v}_k}
\rho e,
\end{equation}
where $e$ is the specific internal energy of the medium.
One supposes that the components have the same velocity $u$ and that no mass transfer occurs between the species. If one notes
$\rho \bW = [\rho Y, \rho u, \rho (e + u^2/2) ]^T$,
$\bT(\bW) 
= [0, 0 , P, P u]^T$
 then the flow is governed by 
\begin{equation}
\partial_t\rho + \partial_x(\rho u) = 0
,\qquad
\partial_t (\rho\bW) + \partial_x (\rho\bW u + \bT) =0.
\label{eq: two fluid model}
\end{equation}
System~\eqref{eq: two fluid model} is hyperbolic provided that $\gamma_k>1$ and it is equipped with jump relations that enable the definition of weak solutions that verifies the transport equation
\begin{equation}
\partial_t Y + u \partial_x Y = 0.
\label{eq: transport of Y}
\end{equation}
Although this model is equipped with a mixture model (based on oversimplified assumptions), if one chooses an initial condition such that $Y(x,t=0)\in\{0,1\}$, then no physical mixing should occur in the domain as \eqref{eq: transport of Y} guarantees that $Y(x,t)\in\{0,1\}$ for $t>0$. In this sense, \eqref{eq: two fluid model} can be used as a model for a flows involving two compressible fluids separated by a sharp interface.

The decoupling between transport and other phenomenon can be achieved thanks to a Lagrange-Remap method. Let us note $\xi$ the Lagrangian space coordinate defined by $\text{d}\xi(t;x_0)/\text{d}t = u(\xi(t;x_0),t)$ with $\xi(t=0;x_0) = x_0$. If $(x,t)\mapsto b$ is any fluid parameter, we note $(\xi,t)\mapsto b^\Lag$ the Lagrangian field associated with $b$ by 
$b^\Lag(\xi(t;x),t) = b(x,t)$. System~\eqref{eq: two fluid model} can be expressed in the so-called Lagrangian reference frame as follows
\begin{equation}
\rho^\Lag(\xi,0) \partial_t( 1/ \rho^\Lag) - \partial_\xi u^\Lag = 0
,\qquad
\rho^\Lag(\xi,0)\partial_t \bW^\Lag + \partial_\xi \bT^\Lag =0.
\label{eq: two fluid model lag}
\end{equation}
Given a set of discrete values $(\rho,\rho\bW)_j^n$ that represent an approximation of the fluid state at instant $t=t^n$ within the cell $i$, the Lagrange-Remap method is a two-step algorithm \cite{godrav1996,Despres2010-book}: first, we update the discrete unknowns to a value
$(\rho,\rho\bW)_j^\Lag$ by approximating the solution of
\eqref{eq: two fluid model lag} over $[t^n,t^n+\Delta t]$.
Let us remark that the evolution equation for $Y$ expressed in \eqref{eq: two fluid model lag} boils down to $\partial_t Y = 0$, therefore it is 
reasonable to expect that $Y^\Lag_i = Y^n_i$.
 The second step
updates the fluid parameter to their values $(\rho,\rho\bW)_j^{n+1}$ by remapping the Lagrange values $(\rho,\rho\bW)_j^\Lag$ onto the Eulerian mesh. It can read as follows
\begin{subequations}
\begin{align}
\rho^{n+1}_j  - \rho^{n}_j 
+ \frac{\Delta t}{\Delta x}
(\rho_{j+1/2}^\Lag u^{n}_{j+1/2} - \rho_{j-1/2}^\Lag u^{n}_{j-1/2}) =0
\label{eq: rho remap step}
,\\
(\rho\bW)^{n+1}_j  - \rho^{n}_j \bW^\Lag_j
+ \frac{\Delta t}{\Delta x}
((\rho\bW)_{j+1/2}^\Lag u^{n}_{j+1/2} 
- 
(\rho\bW)_{j-1/2}^\Lag  u^{n}_{j-1/2}) =0.
\label{eq: W remap step}
\end{align}
\label{eq: remap}
\end{subequations}
The values $u^{n}_{j-1/2}$ are approximations of the material velocity of the fluid at the cell interface $x=x_{j+1/2}$ that can be estimated with the discretization of \eqref{eq: two fluid model lag}. One can therefore consider that $u^{n}_{j-1/2}$ is known when performing \eqref{eq: remap}. The only missing ingredient for achieving the remap procedure is thus the definition of the variable flux 
$(\rho,\rho\bW)_{j+1/2}^\Lag = [\rho,\rho Y, \rho u, \rho (e+u^2/2)]_{j+1/2}^\Lag$. For this problem it is clear that the anti-diffusive mechanism should concern the variable $Y$ whose discontinuity carries the location of material interface between the fluids. Suppose given a definition for the fluxes $b_{j+1/2}^\Lag$ that is consistent for $b\in\{\rho, \rho u, \rho (e+u^2/2)\}$ and that $\rho_j^n>0$ and $\rho_{j+1/2}^\Lag>0$.
Let us note $\mm_{j+1/2}^n= \min(Y^n_j,Y^n_{j+1})$, $\MM_{j+1/2}^n= \max(Y^n_j,Y^n_{j+1})$. Following the ideas introduced in section~\ref{section : inequality and anti-diffusion} in the case of pure transport problem, we aim at defining a flux $Y_{j+1/2}^\Lag$ that fulfills two requirements. 
\begin{itemize}
\item[(i)] $Y_{j+1/2}^\Lag \in [\mm_{j+1/2}^n, \MM_{j+1/2}^n]$;
\item[(ii)] the choice of $Y_{j+1/2}^\Lag$ and \eqref{eq: W remap step} should ensure a discrete maximum principle for $Y$ in the cell $i$ (resp. $i+1$) if $u^n_{j+1/2}>0$ and $u^n_{j-1/2}>0$ 
(resp. $u^n_{j+1/2}<0$ and $u^n_{j+3/2}<0$). 
\end{itemize}
For the sake of simplicity, we suppose that $\rho_{j+1/2}^\Lag$ is defined by the upwind flux, \textit{i.e.}
$
\rho_{j+1/2}^\Lag u^n_{j+1/2} 
= \rho_{j}^\Lag (u^n_{j+1/2})^+ + \rho_{j+1}^\Lag (u^n_{j+1/2})^-$,
then we can define the real interval $[d_{j+1/2},D_{j+1/2}]$ as follows.
\begin{itemize}
\item If  $u^n_{j+1/2}>0$ and $u^n_{j+1/2}>0$ (resp. $u^n_{j+1/2}<0$), 
we set
\begin{align*}
d_{j+1/2} &= 
Y_j^n + (\MM_{j-1/2} - Y^n_j)
\left[ 1 - \frac{\Delta x}{u_{j+1/2}^n\Delta t }\right]
 \quad
 \text{(resp. $d_{j+1/2} =  Y_j^n$ ),}
\\
D_{j+1/2} &= 
Y_j^n + (\mm_{j-1/2} - Y^n_j)
\left[ 1 - \frac{\Delta x}{u_{j+1/2}^n\Delta t }\right]
 \quad
 \text{(resp. $D_{j+1/2} =  Y_j^n$ ).}
\end{align*}

\item If  $u^n_{j+1/2}<0$ and $u^n_{j+3/2}<0$ (resp. $u^n_{j+3/2}>0$), 
we set
\begin{align*}
d_{j+1/2} &= 
Y_{j+1}^n + (\MM_{j+3/2} - Y^n_{j+1})
\left[ 1 + \frac{\Delta x}{u_{j+1/2}^n\Delta t }\right]
 \quad
 \text{(resp. $d_{j+1/2} =  Y_{j+1}^n$ ),}
\\
D_{j+1/2} &= 
Y_{j+1}^n + (\mm_{j+3/2} - Y^n_{j+1})
\left[ 1 - \frac{\Delta x}{u_{j+1/2}^n\Delta t }\right]
 \quad
 \text{(resp. $D_{j+1/2} =  Y_{j+1}^n$ ).}
\end{align*}
\end{itemize}
Let us note $ [\omega_{j+1/2}, \Omega_{j+1/2}] = [\mm_{j+1/2}^n, \MM_{j+1/2}^n] \cap [d_{j+1/2}, D_{j+1/2}]$.
Under the CFL condition
\begin{equation}
|u_{j+1/2}^n|\Delta t / \Delta x < 1,
\label{eq: CFL remap}
\end{equation}
one can check that $[\omega_{j+1/2}, \Omega_{j+1/2}] \neq \emptyset$ as 
$Y_{j}^n$ (resp. $Y_{j+1}^n$) belongs to  $[\omega_{j+1/2}, \Omega_{j+1/2}]$ 
if  $u^n_{j+1/2}>0$  (resp.  $u^n_{j+1/2}<0$ ). Choosing 
$Y_{j+1/2}^\Lag \in [\omega_{j+1/2}, \Omega_{j+1/2}]$ ensures that (i) and (ii) are verified under the condition~\eqref{eq: CFL remap}. In order to enable a sharp transport of $Y$, one just need to use the limited downwind choice within the interval $[\omega_{j+1/2}, \Omega_{j+1/2}]$, which boils down to set
\begin{equation}
Y_{j+1/2}^\Lag = \min(\max(\omega_{j+1/2},Y_\text{down}^\Lag),\Omega_{j+1/2}) 
,
\label{eq: limited downwind choice two-phase model}
\end{equation}
where $Y_\text{down}^\Lag = Y_{j+1}^\Lag$ 
(resp. $Y_\text{down}^\Lag = Y_{j}^\Lag$)
if $u^n_{j+1/2}>0$ (resp. $u^n_{j+1/2}<0$).

A numerical scheme based on a Finite Volume approximation of \eqref{eq: two fluid model lag} and \eqref{eq: remap} with the limited downwind choice 
\eqref{eq: limited downwind choice two-phase model} was studied in \cite{lagoutiere2000} for the model described in this section. It is worth mentioning that up to a careful discretization choice for \eqref{eq: two fluid model lag} 
the overall algorithm is conservative with respect to $(\rho,\rho\bW)$. Let us also emphasize that the algorithm presented in this section is difficult to use 
in practice: spurious pressure and velocity oscillations at the material interface may occurs, which is a common issue for this type of problems \cite{abgrall1996}.
The same method was applied to similar two-phase models with an alternate mixture law in \cite{lagoutiere2000,DesLag2007} that guarantees that constant pressure and velocity profiles are preserved.

\subsection{Example of other evolution equation involving sharp interfaces}
It is not possible to give an exhaustive list of all possible sharpening techniques implementation, we will try to give hereafter an overview of the works that have been achieved the past years that is inevitably incomplete.

The approach of section~\ref{section: system with dalton law} has been successfully extended to other systems like the five-equation model of \cite{massoni,AkClKo2002} in \cite{KoLag2010} and also for compressibles flows involving an arbitrary number of components separated by interfaces~\cite{jaouen-lagoutiere,billaud-kokh}. Other techniques may be used to sharpen front in systems with interface. For example, considering again system~\eqref{eq: two fluid model}, one can discretize directly the transport equation~\eqref{eq: transport of Y} with the limited downwind scheme of section~\ref{section : inequality and anti-diffusion} and use a classical Finite-Volume discretization for $\rho$, $\rho u$ and $\rho (e+u^2/2)$, at the cost of deriving a non-conservative numerical scheme. Other sharpening techniques can also be used for compressible two-phase flows with interface similar to \eqref{eq: two fluid model}: the THINC method that was first developed for incompressible flows~\cite{XiaoHonmaKono2005} has been adapted in \cite{shyue-xiao} to the five-equation model studied in \cite{AkClKo2002}. This method relies on controlling the spreading of the material interface thanks to an hyperbolic tangent profile. As mentioned in section~\ref{section: Glimm}, Glimm's method has also been used for discretizing sharply the evolution of an interface. Indeed, it is possible to sharply let evolve contact discontinuities in a system by 
providing a dedicated treatment based on a Glimm type random choice method \cite{chalons2007,ChalonsGoatin2008}. In \cite{Bachmann2013275} a random choice method is within a Lagrange-Remap strategy to perform the Remap step while preserving sharp profiles.
The limited downwind strategy has been implemented to describe interface that are not solely passively advected like problems of reacting gas flows\cite{Tang2014}.
A VOF-type reconstruction that relies on a level set description of the interface is proposed in \cite{hu2006} for the simulation of two-component compressible flows.

\subsection{Cut-cells and CFL condition}

Taking as a principle that sharpening techniques have the ability to reconstruct 
interfaces, it appears that an interface which moves dynamically in a Cartesian mesh 
may cut cells into smaller cells. Of course it is most of the time 
only a geometrical interpretation. However it has the unfortunate consequence
that these small cut cells may have a dramatic influence on the 
CFL conditions through a complex nonlinear interaction of the parts of the global algorithm (note nevertheless that it is not the case with the Vofire algorithm).
This feature is difficult  to analyze rigorously in the context of  sharpening methods. In practice 
one observes a posteriori the stability or the instability of the scheme/code.
We refer to the chapter \cite{berger} in this volume for a comprehensive presentation of the topic.

\bibliographystyle{alpha}
\bibliography{paper-review-sharpening-interface}

\newcommand{\etalchar}[1]{$^{#1}$}
\begin{thebibliography}{HKAH06}

\bibitem[Abg88]{abgrall2}
R.~Abgrall.
\newblock Generalization of the roe scheme for the computation of mixture of
  perfect gases.
\newblock {\em Rech. A\'erospatiale (English edition)}, 6:31--43, 1988.

\bibitem[Abg96]{abgrall1996}
R.~Abgrall.
\newblock How to prevent pressure oscillations in multicomponent flow
  calculations: a quasi-conservative approach.
\newblock {\em J. Comput. Phys.}, 125(1):150--160, 1996.

\bibitem[AC16]{aguillon-chalons}
N.~Aguillon and C.~Chalons.
\newblock Nondiffusive conservative schemes based on approximate riemann
  solvers for {L}agrangian gas dynamics.
\newblock {\em ESAIM: Mathematical Modelling and Numerical Analysis}, to
  appear, 2016.

\bibitem[ACK02]{AkClKo2002}
G.~Allaire, S.~Clerc, and S.~Kokh.
\newblock A five-equation model for the simulation of interfaces between
  compressible fluids.
\newblock {\em J. Comput. Phys.}, 181(2):577--616, 2002.

\bibitem[Agu16]{aguillon}
N.~Aguillon.
\newblock Capturing nonclassical shocks in nonlinear elastodynamic with a
  conservative finite volume scheme.
\newblock {\em Interfaces Free Bound.}, to appear, 2016.

\bibitem[BCLL08]{boutin-chalons-lagoutiere-lefloch}
B.~Boutin, C.~Chalons, F.~Lagouti{\`e}re, and P.~G. LeFloch.
\newblock Convergent and conservative schemes for nonclassical solutions based
  on kinetic relations. {I}.
\newblock {\em Interfaces Free Bound.}, 10(3):399--421, 2008.

\bibitem[Ber84]{berger}
M.~Berger.
\newblock In {\em Computing methods in applied sciences and engineering, {VI}
  ({V}ersailles, 1983)}, pages 491--492. North-Holland, Amsterdam, 1984.

\bibitem[BFK14]{billaud-kokh}
M.~Billaud~Friess and S.~Kokh.
\newblock Simulation of sharp interface multi-material flows involving an
  arbitrary number of components through an extended five-equation model.
\newblock {\em J. Comput. Phys.}, 273:488--519, 2014.

\bibitem[BHJ{\etalchar{+}}13]{Bachmann2013275}
M.~Bachmann, P.~Helluy, J.~Jung, H.~Mathis, and S.~M{\"u}ller.
\newblock Random sampling remap for compressible two-phase flows.
\newblock {\em Comput. \& Fluids}, 86:275--283, 2013.

\bibitem[BTVG10]{tran1}
J.~Bohbot, Q.~H. Tran, A.~Velghe, and N.~Gillet.
\newblock A multi-dimensional spatial scheme for massively parallel
  compressible turbulent combustion simulation.
\newblock In {\em Proceedings of the V European Conference on Computational
  Fluid Dynamics, ECCOMAS CFD 2010}, pages 1--20. J. C. F. Pereira, A. Sequeira
  and J. M. C. Pereira (Eds), 2010.

\bibitem[CC12]{chalons-coquel}
C.~Chalons and F.~Coquel.
\newblock Computing material fronts with a {L}agrange-projection approach.
\newblock In {\em Hyperbolic problems---theory, numerics and applications.
  {V}olume 1}, volume~17 of {\em Ser. Contemp. Appl. Math. CAM}, pages
  346--356. World Sci. Publishing, Singapore, 2012.

\bibitem[CG08]{ChalonsGoatin2008}
C.~Chalons and P.~Goatin.
\newblock Godunov scheme and sampling technique for computing phase transitions
  in traffic flow modeling.
\newblock {\em Interfaces Free Bound.}, 10(2):197--221, 2008.

\bibitem[Cha07]{chalons2007}
C.~Chalons.
\newblock Numerical approximation of a macroscopic model of pedestrian flows.
\newblock {\em SIAM J. Sci. Comput.}, 29(2):539--555 (electronic), 2007.

\bibitem[Cho76]{chorin}
A.~J. Chorin.
\newblock Random choice solution of hyperbolic systems.
\newblock {\em J. Computational Phys.}, 22(4):517--533, 1976.

\bibitem[CL03]{chalons-lefloch}
C.~Chalons and P.~G. LeFloch.
\newblock Computing undercompressive waves with the random choice scheme.
  {N}onclassical shock waves.
\newblock {\em Interfaces Free Bound.}, 5(2):129--158, 2003.

\bibitem[CM11]{MR2787572}
R.~Chen and D.-K. Mao.
\newblock Entropy-{TVD} scheme for nonlinear scalar conservation laws.
\newblock {\em J. Sci. Comput.}, 47(2):150--169, 2011.

\bibitem[Col82]{colella}
P.~Colella.
\newblock Glimm's method for gas dynamics.
\newblock {\em SIAM J. Sci. Statist. Comput.}, 3(1):76--110, 1982.

\bibitem[CPT12]{MR2877636}
Y.~G. Chen, W.~G. Price, and P.~Temarel.
\newblock An anti-diffusive volume of fluid method for interfacial fluid flows.
\newblock {\em Internat. J. Numer. Methods Fluids}, 68(3):341--359, 2012.

\bibitem[Des08]{despresjll}
B.~Despr{\'e}s.
\newblock Finite volume transport schemes.
\newblock {\em Numer. Math.}, 108(4):529--556, 2008.

\bibitem[Des09]{desstab}
B.~Despr{\'e}s.
\newblock Uniform asymptotic stability of {S}trang's explicit compact schemes
  for linear advection.
\newblock {\em SIAM J. Numer. Anal.}, 47(5):3956--3976, 2009.

\bibitem[Des10]{Despres2010-book}
B.~Despr{\'e}s.
\newblock {\em Lois de conservations eul\'eriennes, lagrangiennes et m\'ethodes
  num\'eriques}, volume~68 of {\em Math\'ematiques \& Applications (Berlin)
  [Mathematics \& Applications]}.
\newblock Springer-Verlag, Berlin, 2010.

\bibitem[DL01a]{DesLag02}
B.~Despr{\'e}s and F.~Lagouti{\`e}re.
\newblock Contact discontinuity capturing schemes for linear advection and
  compressible gas dynamics.
\newblock {\em J. Sci. Comput.}, 16(4):479--524 (2002), 2001.

\bibitem[DL01b]{lvd}
B.~Despr{\'e}s and F.~Lagouti{\`e}re.
\newblock Generalized {H}arten formalism and longitudinal variation diminishing
  schemes for linear advection on arbitrary grids.
\newblock {\em M2AN Math. Model. Numer. Anal.}, 35(6):1159--1183, 2001.

\bibitem[DL07]{DesLag2007}
B.~Despr{\'e}s and F.~Lagouti{\`e}re.
\newblock Numerical resolution of a two-component compressible fluid model with
  interfaces.
\newblock {\em Prog. Comput. Fluid Dyn.}, 7(6):295--310, 2007.

\bibitem[DL11]{delarue-lagoutiere}
F.~Delarue and F.~Lagouti{\`e}re.
\newblock Probabilistic analysis of the upwind scheme for transport equations.
\newblock {\em Arch. Ration. Mech. Anal.}, 199(1):229--268, 2011.

\bibitem[DLLM10]{vofire}
B.~Despr{\'e}s, F.~Lagouti{\`e}re, E.~Labourasse, and I.~Marmajou.
\newblock An antidissipative transport scheme on unstructured meshes for
  multicomponent flows.
\newblock {\em Int. J. Finite Vol.}, 7(1):36, 2010.

\bibitem[DM96]{dubois}
F.~Dubois and G.~Mehlman.
\newblock A non-parameterized entropy correction for {R}oe's approximate
  {R}iemann solver.
\newblock {\em Numer. Math.}, 73(2):169--208, 1996.

\bibitem[GL85]{good-lev}
J.~B. Goodman and R.~J. LeVeque.
\newblock On the accuracy of stable schemes for {$2$}{D} scalar conservation
  laws.
\newblock {\em Math. Comp.}, 45(171):15--21, 1985.

\bibitem[Gli65]{glimm}
J.~Glimm.
\newblock Solutions in the large for nonlinear hyperbolic systems of equations.
\newblock {\em Comm. Pure Appl. Math.}, 18:697--715, 1965.

\bibitem[GLT13]{goudon-lagoutiere-tine}
T.~Goudon, F.~Lagouti{\`e}re, and L.~M. Tine.
\newblock Simulations of the {L}ifshitz-{S}lyozov equations: the role of
  coagulation terms in the asymptotic behavior.
\newblock {\em Math. Models Methods Appl. Sci.}, 23(7):1177--1215, 2013.

\bibitem[GR96]{godrav1996}
E.~Godlewski and P.-A. Raviart.
\newblock {\em Numerical approximation of hyperbolic systems of conservation
  laws}, volume 118 of {\em Applied Mathematical Sciences}.
\newblock Springer-Verlag, New York, 1996.

\bibitem[Har84]{harten1}
A.~Harten.
\newblock On a class of high resolution total-variation-stable
  finite-difference schemes.
\newblock {\em SIAM J. Numer. Anal.}, 21(1):1--23, 1984.
\newblock With an appendix by Peter D. Lax.

\bibitem[HJ13]{helluy-jung2}
P.~Helluy and J.~Jung.
\newblock Opencl simulations of two-fluid compressible flows with a random
  choice method.
\newblock {\em IJFV}, 10:1–38, 2013.

\bibitem[HJ14]{helluy-jung1}
P.~Helluy and J.~Jung.
\newblock Two-fluid compressible simulations on {GPU} cluster.
\newblock In {\em Congr\`es {SMAI} 2013}, volume~45 of {\em ESAIM Proc.
  Surveys}, pages 349--358. EDP Sci., Les Ulis, 2014.

\bibitem[HKAH06]{hu2006}
X.~Y. Hu, B.~C. Khoo, N.~A. Adams, and F.-L. Huang.
\newblock A conservative interface method for compressible flows.
\newblock {\em J. Comput. Phys.}, 219(2):553--578, 2006.

\bibitem[HL94]{hou-lefloch}
T.~Y. Hou and P.~G. LeFloch.
\newblock Why nonconservative schemes converge to wrong solutions: error
  analysis.
\newblock {\em Math. Comp.}, 62(206):497--530, 1994.

\bibitem[HN81]{vof}
C.-W. Hirt and B.~D. Nichols.
\newblock Volume of fluid (vof) method for the dynamics of free boundaries.
\newblock {\em J. Comput. Phys.}, 39:201--225, 1981.

\bibitem[IS83]{iserles1}
A.~Iserles and G.~Strang.
\newblock The optimal accuracy of difference schemes.
\newblock {\em Trans. Amer. Math. Soc.}, 277(2):779--803, 1983.

\bibitem[JL07]{jaouen-lagoutiere}
S.~Jaouen and F.~Lagouti{\`e}re.
\newblock Numerical transport of an arbitrary number of components.
\newblock {\em Comput. Methods Appl. Mech. Engrg.}, 196(33-34):3127--3140,
  2007.

\bibitem[KL10]{KoLag2010}
S.~Kokh and F.~Lagouti{\`e}re.
\newblock An anti-diffusive numerical scheme for the simulation of interfaces
  between compressible fluids by means of a five-equation model.
\newblock {\em J. Comput. Phys.}, 229(8):2773--2809, 2010.

\bibitem[KM05a]{kuz-tur3}
D.~Kuzmin and M.~M{\"o}ller.
\newblock Algebraic flux correction. {I}. {S}calar conservation laws.
\newblock In {\em Flux-corrected transport}, Sci. Comput., pages 155--206.
  Springer, Berlin, 2005.

\bibitem[KM05b]{kuz-tur2}
D.~Kuzmin and M.~M{\"o}ller.
\newblock Algebraic flux correction. {II}. {C}ompressible {E}uler equations.
\newblock In {\em Flux-corrected transport}, Sci. Comput., pages 207--250.
  Springer, Berlin, 2005.

\bibitem[Lag00]{lagoutiere2000}
F.~Lagouti\`ere.
\newblock {\em Mod\'elisation math\'ematique et r\'esolution num\'erique de
  probl\`emes de fluides compressibles \`a plusieurs constituants}.
\newblock PhD thesis, Universit\'e Paris VI, 2000.

\bibitem[LeF02]{lefloch}
P.~G. LeFloch.
\newblock {\em Hyperbolic systems of conservation laws}.
\newblock Lectures in Mathematics ETH Z\"urich. Birkh\"auser Verlag, Basel,
  2002.
\newblock The theory of classical and nonclassical shock waves.

\bibitem[LF89]{larrouturou1989}
B.~Larrouturou and L.~F{\'e}zoui.
\newblock On the equations of multi-component perfect or real gas inviscid
  flow.
\newblock In {\em Nonlinear hyperbolic problems ({B}ordeaux, 1988)}, volume
  1402 of {\em Lecture Notes in Math.}, pages 69--98. Springer, Berlin, 1989.

\bibitem[lR77]{Leroux}
A.-Y. le~Roux.
\newblock A numerical conception of entropy for quasi-linear equations.
\newblock {\em Math. Comp.}, 31(140):848--872, 1977.

\bibitem[LW60]{LaxWen60}
P.~Lax and B.~Wendroff.
\newblock Systems of conservation laws.
\newblock {\em Comm. Pure Appl. Math.}, 13:217--237, 1960.

\bibitem[MSNA02]{massoni}
J.~Massoni, R.~Saurel, B.~Nkonga, and R.~Abgrall.
\newblock Some models and {E}ulerian methods for interface problems between
  compressible fluids with heat transfer.
\newblock {\em Int. J. Heat Mass Transfer}, 45(6):1287--1307, 2002.

\bibitem[MTF10]{tran2}
A.~Michel, Q.~H. Tran, and G.~Favennec.
\newblock A genuinely one-dimensional upwind scheme with accuracy enhancement
  for multidimensional advection problems.
\newblock In {\em ECMOR XII – 12 th European Conference on the Mathematics of
  Oil Recovery 6-9 September 2010, Oxford, UK}, pages 1--21. European
  Association of Geoscientists \& Engineers, 2010.

\bibitem[NW76]{Slic}
W.~F. Noh and P.~Woodward.
\newblock Slic (simple line interface calculation).
\newblock {\em Communications in Mathematical Sciences}, 59:57--70, 1976.

\bibitem[OF03]{osher}
S.~Osher and R.~Fedkiw.
\newblock {\em Level set methods and dynamic implicit surfaces}, volume 153 of
  {\em Applied Mathematical Sciences}.
\newblock Springer-Verlag, New York, 2003.

\bibitem[OS88]{osher-sethian}
S.~Osher and J.~A. Sethian.
\newblock Fronts propagating with curvature-dependent speed: algorithms based
  on {H}amilton-{J}acobi formulations.
\newblock {\em J. Comput. Phys.}, 79(1):12--49, 1988.

\bibitem[Roe85]{roe0}
P.~L. Roe.
\newblock Some contributions to the modelling of discontinuous flows.
\newblock In {\em Large-scale computations in fluid mechanics, {P}art 2 ({L}a
  {J}olla, {C}alif., 1983)}, volume~22 of {\em Lectures in Appl. Math.}, pages
  163--193. Amer. Math. Soc., Providence, RI, 1985.

\bibitem[RTT08]{MR2549232}
G.~Russo, E.~F. Toro, and V.~A. Titarev.
\newblock A{DER}-{R}unge-{K}utta schemes for conservation laws in one space
  dimension.
\newblock In {\em Hyperbolic problems: theory, numerics, applications}, pages
  929--936. Springer, Berlin, 2008.

\bibitem[Set96]{sethian}
J.~A. Sethian.
\newblock {\em Level set methods}, volume~3 of {\em Cambridge Monographs on
  Applied and Computational Mathematics}.
\newblock Cambridge University Press, Cambridge, 1996.
\newblock Evolving interfaces in geometry, fluid mechanics, computer vision,
  and materials science.

\bibitem[Shu09]{MR2481112}
C.-W. Shu.
\newblock High order weighted essentially nonoscillatory schemes for convection
  dominated problems.
\newblock {\em SIAM Rev.}, 51(1):82--126, 2009.

\bibitem[SST15]{Schmidtmann}
B.~Schmidtmann, B.~Seibold, and M.~Torrilhon.
\newblock Relations between {WENO3} and third-order limiting in finite volume
  methods.
\newblock {\em Journal of Scientific Computing}, pages 1--29, 2015.

\bibitem[Str68]{strang}
G.~Strang.
\newblock On the construction and comparison of difference schemes.
\newblock {\em SIAM J. Numer. Anal.}, 5:506--517, 1968.

\bibitem[Swe84]{sweb1}
P.~K. Sweby.
\newblock High resolution schemes using flux limiters for hyperbolic
  conservation laws.
\newblock {\em SIAM J. Numer. Anal.}, 21(5):995--1011, 1984.

\bibitem[Swe85]{sweb2}
P.~K. Sweby.
\newblock High resolution {TVD} schemes using flux limiters.
\newblock In {\em Large-scale computations in fluid mechanics, {P}art 2 ({L}a
  {J}olla, {C}alif., 1983)}, volume~22 of {\em Lectures in Appl. Math.}, pages
  289--309. Amer. Math. Soc., Providence, RI, 1985.

\bibitem[SX14]{shyue-xiao}
K.-M. Shyue and F.~Xiao.
\newblock An {E}ulerian interface sharpening algorithm for compressible
  two-phase flow: the algebraic {THINC} approach.
\newblock {\em J. Comput. Phys.}, 268:326--354, 2014.

\bibitem[TBC14]{Tang2014}
K.~Tang, A.~Beccantini, and C.~Corre.
\newblock Combining discrete equations method and upwind downwind-controlled
  splitting for non-reacting and reacting two-fluid computations: two
  dimensional case.
\newblock {\em Comput. \& Fluids}, 103:132--155, 2014.

\bibitem[TK05]{kuz-tur}
S.~Turek and D.~Kuzmin.
\newblock Algebraic flux correction. {III}. {I}ncompressible flow problems.
\newblock In {\em Flux-corrected transport}, Sci. Comput., pages 251--296.
  Springer, Berlin, 2005.

\bibitem[Tor97]{Tor97}
E.~F. Toro.
\newblock {\em Riemann solvers and numerical methods for fluid dynamics}.
\newblock Springer-Verlag, Berlin, 1997.
\newblock A practical introduction.

\bibitem[TT05]{MR2219366}
E.~F. Toro and V.~A. Titarev.
\newblock T{VD} fluxes for the high-order {ADER} schemes.
\newblock {\em J. Sci. Comput.}, 24(3):285--309, 2005.

\bibitem[TT07]{MR2337583}
V.~A. Titarev and E.~F. Toro.
\newblock Analysis of {ADER} and {ADER}-{WAF} schemes.
\newblock {\em IMA J. Numer. Anal.}, 27(3):616--630, 2007.

\bibitem[WB76]{beam:warming}
R.~F. Warming and R.~M. Beam.
\newblock Upwind second-order difference schemes and applications in
  aerodynamic flows.
\newblock {\em AIAA J.}, 14(9):1241--1249, 1976.

\bibitem[XHK05]{XiaoHonmaKono2005}
F.~Xiao, Y.~Honma, and T.~Kono.
\newblock A simple algebraic interface capturing scheme using hyperbolic
  tangent function.
\newblock {\em International Journal for Numerical Methods in Fluids},
  48(9):1023--1040, 2005.

\bibitem[XS05]{cws2}
Z.~Xu and C.-W. Shu.
\newblock Anti-diffusive high order {WENO} schemes for {H}amilton-{J}acobi
  equations.
\newblock {\em Methods Appl. Anal.}, 12(2):169--190, 2005.

\bibitem[XS06]{cws1}
Z.~Xu and C.-W. Shu.
\newblock Anti-diffusive finite difference {WENO} methods for shallow water
  with transport of pollutant.
\newblock {\em J. Comput. Math.}, 24(3):239--251, 2006.

\bibitem[You84]{youngs}
D.~L. Youngs.
\newblock An interface tracking method for a 3{D} {E}ulerian hydrodynamics
  code.
\newblock Technical Report Technical Report 44/92/35, AWRE, 1984.

\bibitem[YWH84]{harten2}
H.~C. Yee, R.~F. Warming, and A.~Harten.
\newblock On a class of {TVD} schemes for gas dynamic calculations.
\newblock In {\em Computing methods in applied sciences and engineering, {VI}
  ({V}ersailles, 1983)}, pages 491--492. North-Holland, Amsterdam, 1984.

\end{thebibliography}

\end{document}